\documentclass{amsart}
\usepackage{amsmath}
\usepackage[dvips]{graphicx}

\newtheorem{theorem}{Theorem}
\newtheorem{lemma}[theorem]{Lemma}

\newtheorem{corollary}[theorem]{Corollary}

\theoremstyle{definition}
\newtheorem{definition}[theorem]{Definition}
\newtheorem{remark}[theorem]{Remark}

\newcommand{\cpk}{\mathbb{CP}^2}
\newcommand{\cpkk}{\overline{\mathbb{CP}}^2}
\newcommand{\zz}{\mathbb{Z}}

\newcommand{\rr}{\mathbb{R}}

\begin{document}

\title{Non-isotopic Symplectic Tori in the Same Homology Class}

\author{Tolga Etg\"u}
\address{Department of Mathematics and Statistics, McMaster University,
Hamilton, Ontario L8S 4K1, Canada}
\email{etgut@math.mcmaster.ca}

\author{B. Doug Park}
\address{Department of Pure Mathematics,
University of Waterloo, Waterloo, Ontario N2L 3G1, Canada}
\email{bdpark@math.uwaterloo.ca}
%\curraddr{}
%\thanks{}


\date{October 28, 2002.  Revised on December 12, 2002}

\begin{abstract}
For any pair of integers\/ $n\geq 1$\/ and\/ $q\geq 2$, we
construct an infinite family of mutually non-isotopic symplectic
tori representing the homology class\/ $q[F]$\/ of an elliptic
surface $E(n)$, where\/ $[F]$\/ is the homology class of the
fiber.  We also show how such families can be non-isotopically
and symplectically embedded into a more general class of
symplectic $4$-manifolds.
\end{abstract}

\maketitle

\section{Introduction}

In a fixed homology class of a complex surface, there are at
most finitely many complex curves up to smooth isotopy. In
contrast, there exist examples of homology classes in symplectic
$4$-manifolds that are represented by  infinitely many
non-isotopic connected symplectic surfaces. In
\cite{fs:non-isotopic} Fintushel and Stern proved that
$q[F]\in H_2(X; \zz)$\/ is such a class provided that $q \geq 4$
is an even integer and $F$\/ is a symplectic $c$-embedded torus, i.e.
a homologically nontrivial torus of self-intersection $0$ with
first homology generated
by two circles which bound disks of self-intersection $-1$ in $X$,
in a simply-connected symplectic $4$-manifold $X$
(e.g. we can take $X$\/ to be an elliptic surface $E(n)$
and $F$ to be a regular fiber).
In this paper we generalize their result by constructing
infinitely many non-isotopic symplectic tori in
$q[F] \in H_2(X; \zz)$\/ for each $q \geq 3$ (even or odd) in
case when $F$\/ is an essentially embedded
(cf.$\;$Definition~\ref{def:essential})
symplectic torus in a
symplectic $4$-manifold $X$\/ satisfying\/ $[F]^2 =0$\/ and\/
$H^1(X \setminus\nu(F); \zz) = 0$. Moreover, as a consequence
of our calculation, it will be apparent that no two of these
tori are equivalent under the action of {\it Diff}$\,(X)$, the
group of self-diffeomorphisms of $X$.
In some special cases our result can be slightly improved
to include $2[F]$ as in the following theorem.

\begin{theorem}\label{theorem:1}
For any pair of integers\/ $n\geq 1$\/ and\/ $q\geq 2$, there
exists an infinite family of mutually non-isotopic symplectic tori
representing the homology class\/ $q[F]$\/ of an elliptic surface
$E(n)$, where\/ $[F]$\/ is the homology class of the fiber.
\end{theorem}

It should be noted that in \cite{vidussi:non-isotopic}, Vidussi
has constructed such families for every positive multiple of the
fiber class in\/ $E(n)$\/ provided that\/ $n\geq3$. However, it
still remains an interesting open problem to construct infinite
families of non-isotopic symplectic tori in the fiber class\/ $[F]$\/
of\/ $E(1)$\/ and\/ $E(2)$.  We should also note that our construction 
can be carried out for a more general class of 4-manifolds than those in 
\cite{vidussi:non-isotopic}.

In all of the above examples, families of braids are used to
construct symplectic tori and Seiberg-Witten theory is used to
distinguish them. The family of braids we use in this paper are,
in some sense, the simplest possible and this allowed us not only
to calculate the Seiberg-Witten invariants of the corresponding
$4$-manifolds completely, but also to reprove and extend the
results in \cite{fs:non-isotopic} and \cite{vidussi:non-isotopic}
(except for the case of fiber class itself in\/ $E(n)$\/ for
$n\geq3$).  From our construction, it will  be apparent
that any other collection of reasonably complicated braids will
also give rise to many more new examples of non-isotopic symplectic tori.

In the next section we construct a family of tori,
$\{T_{p,q}\}_{p , q\geq2}$\hspace{1pt},  in $T^2 \times D^2$ which
could be identified with a regular tubular neighborhood of a torus
of self-intersection $0$ in a
$4$-manifold. Then we show how these tori could be symplectically
embedded in $E(n)$. In Section~\ref{section:sw}
we calculate the Seiberg-Witten
invariants of symplectic $4$-manifolds obtained as the fiber sum
of\/ $E(n)$\/ and\/ $E(r)$\/ along\/ $T_{p,q} \subset E(n)$\/ and a
fiber in\/ $E(r)$. As a consequence of these computations
we prove Theorem~\ref{theorem:1}. In Section~\ref{section:generalization}
we explain how
Theorem~\ref{theorem:1} could be generalized to a larger class
of $4$-manifolds (see Theorem~\ref{theorem:2}).

\section{Braid Construction}

We first construct an infinite family of non-isotopic symplectic
tori\/ $\{T_{p,q}\}_{p,q\geq 2}$\/ representing\/ $q[F]$\/ in
$E(1)$, and show that these remain non-isotopic after
fiber-summing with $E(n-1)$. I.e. we will construct an infinite
family of tori that embed non-isotopically into $E(n)$\/ for any\/
$n\geq 1$.

\begin{figure}[!ht]
\begin{center}\hspace{-1.2cm}
\includegraphics[scale=0.4]{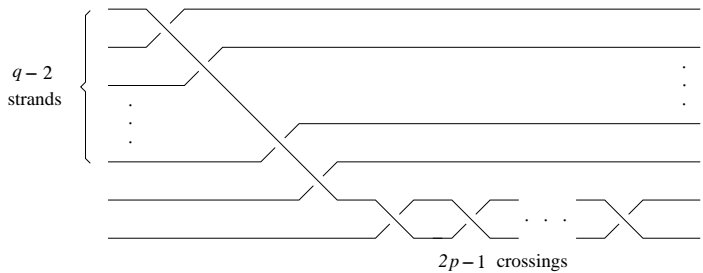}
\end{center}
\caption{Braid $B_{p,q}$ for $p\geq 1$\/ and $q\geq 2$}
\label{fig:braid}
\end{figure}

Even though we are unable to show that\/ $T_{1,q}$\/ is not isotopic to
any torus in $\{ T_{p,q} \:|\; p \geq 2\}$ in general, we will still include
it in our construction and computations for the sake of completeness.
For any pair of integers\/ $p\geq 1$\/
and\/ $q\geq 2$, consider
the $q$-strand braid\/ $B_{p,q}$ in Figure~\ref{fig:braid}.  The
closed braid\/ $\hat{B}_{p,q}$ is the `simplest' knot in\/ $S^3$
with\/ $(2p-1)$ crossings, which is usually denoted by \/ $(2p-1)_1$
in the literature (see e.g. \cite{bz} or \cite{[Ro]}).

\begin{figure}[!ht]
\begin{center}
\includegraphics[scale=0.3]{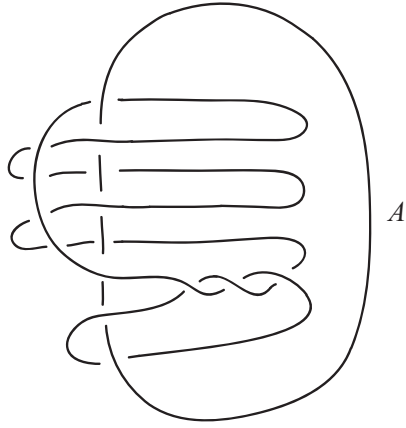}
\end{center}
\caption{Two-component link $L_{2,4}$} \label{fig:L24}
\end{figure}

Let\/ $L_{p,q}$ denote the two-component link in\/ $S^3$ which is
the union of the closed braid $\hat{B}_{p,q}$ and its axis $A$.
See Figure~\ref{fig:L24} for a picture of\/ $L_{2,4}\hspace{1pt}$.

Next we review the generalization of the link surgery construction
of Fintushel and Stern \cite{fs:knots} by Vidussi \cite{vidussi:smooth}.
For an $m$-component link $L$, choose a homology basis of simple
curves $\{\alpha_i, \beta_i \}_{i=1}^{m}$ such that the intersection of
$\alpha_i$ and $\beta_i$ is 1 in the boundary of the
link exterior. We define the link
surgery manifold
\[
E(n_1, \ldots\hspace{-.5pt} ,n_m)_L \: :=\; [\coprod_{i=1}^{m}
E(n_i)\setminus\nu F_i]\hspace{-20pt}\bigcup_{F_i\times\partial
D^2=(S^1\times \alpha_i)\times\beta_i}\hspace{-20pt}
[S^1\times(S^3\setminus \nu L)]\, ,
\]
where $\nu$ denotes the tubular neighborhoods.  Here, the gluing
diffeomorphisms between the boundary 3-tori identify the torus
fiber $F_i$ of $E(n_i)$\/ with\/ $S^1\times \alpha_i\,$, and act
as complex conjugation on the last remaining\/ $S^1$ factor.
Strictly speaking, our gluing construction depends on our choice
of basis\/ $\{\alpha_i, \beta_i \}$, but we will suppress this
dependence in our notation.

As before let $A$\/ denote the unknot which is the axis of the
closed braid $\hat{B}_{p,q}\,$.  Choose the homology basis\/
$(\alpha_1, \beta_1) = (\mu(A),\lambda(A))$\/ for the complement\/
$(S^3 \setminus \nu A)$, where $\mu$ and $\lambda$ denote the
meridian and longitude, respectively.
\begin{lemma}\label{lemma:E(1)}
$E(1)_A = E(1)$.
\end{lemma}

\begin{proof}
Note that the exterior of the unknot\/ $(S^3 \setminus \nu A)$\/
is diffeomorphic to a solid torus\/ $S^1 \times D^2$.  Hence there
is a diffeomorphism between the tubular neighborhood of a regular
fiber $\nu F = T^2 \times D^2$\/ and the Cartesian product\/ $[S^1
\times (S^3 \setminus \nu A)]$.
\end{proof}

Now consider $\hat{B}_{p,q}$ as a knot inside the solid torus
$(S^3 \setminus \nu A)$.  Define $T_{p,q}\, :=S^1\times
\hat{B}_{p,q}\,$, which is a torus embedded in\/ $[S^1 \times (S^3
\setminus \nu A)] \approx \nu F \subset E(1)_A\,$.

\begin{lemma}\label{lemma:symplectic}
The torus\/ $T_{p,q} \subset E(1)$\/ is a symplectic submanifold
and represents the homology class\/ $q[F]$.
\end{lemma}

\begin{proof}
First note that $\hat{B}_{p,q}$ is homologous to\/ $q \mu(A)$\/
in\/ $H_1(S^3\setminus \nu A)$.  Since the fiber $F$\/ is
identified with\/ $S^1\times \mu(A)$\/ in $E(1)_A$\hspace{1pt},
the torus $T_{p,q}$ is homologous to\/ $qF$\/ in\/ $E(1)_A=E(1)$.
To show that $T_{p,q}$ is symplectic, we can proceed as in
\cite{fs:non-isotopic} and \cite{vidussi:non-isotopic}.  The key
step is once again verifying that the knot $\hat{B}_{p,q}\subset
(S^3\setminus \nu A) \approx (S^1\times D^2)$\/ is transverse to
the disks\/ $\{ {\rm pt.} \}\times D^2$.
\end{proof}

Next we show that the torus $T_{p,q}$ can be embedded in $E(n)$
for any $n\geq 1$.  For $n\geq 2$, we can write $E(n)$ as the
fiber sum,
\begin{equation}\label{fiber-sum}
E(n) \,=\: E(1)_A \#_{F'_1=F_2} E(n-1) \: =\; [E(1)_A \setminus \nu
F'_1]\hspace{-20pt} \bigcup_{F'_1\times \partial D^2 = F_2\times
\partial D^2} \hspace{-20pt} [E(n-1)\setminus \nu F_2]\;,
\end{equation}
where $F'_1$ is a regular torus fiber lying outside\/ $[S^1 \times
(S^3 \setminus \nu A)]$.  We immediately see that $T_{p,q}$ can be
embedded into the first fiber summand and hence into\/ $E(n)$.
Without causing too much confusion, we will continue to denote
this embedded torus in $E(n)$\/ by\/ $T_{p,q}\hspace{1pt}$. Note
that Lemma~\ref{lemma:symplectic} continues to hold for the torus\/
$T_{p,q}\subset E(n)$.

Since the torus\/ $T_{p,q}\subset E(n)$\/ has self-intersection
zero, we can form the fiber sum,
\begin{equation}\label{fiber-sum2}
E(n)\#_{T_{p,q}=F} E(r) \: := \; [E(n)\setminus \nu T_{p,q}]
\hspace{-20pt} \bigcup_{T_{p,q}\times\partial D^2=F\times \partial
D^2} \hspace{-20pt} [E(r)\setminus \nu F]\; .
\end{equation}
The gluing diffeomorphism $\varphi$ identifies $T_{p,q}$ with
$F$\/ and acts as complex conjugation on the $\partial D^2$
factor.  Note that when $n\geq 2$\/ 
there is a canonical framing of $T_{p,q}$ in $E(n)$ 
that is inherited from a fixed 
framing of $T_{p,q}$ in $E(1)$ by our fiber sum description 
(\ref{fiber-sum}).  
We remark that the gluing formula in \cite{doug:pft3}
cannot be used to compute the Seiberg-Witten invariant of the
fiber sum (\ref{fiber-sum2}), since the triple\/ $([E(n)\setminus\nu
T_{p,q}],[E(r)\setminus \nu F],\varphi)$\/ is not admissible in
the sense of Definition~8 in \cite{doug:pft3}.  However, the key
observation is that this fiber sum can be expressed as a link
surgery manifold.

\begin{lemma}\label{lemma:main}
$E(n)\#_{T_{p,q}=F} E(r) = \, E(n,r)_{L_{p,q}}\hspace{1pt}$.
\end{lemma}

\begin{proof}
For the construction of\/ $E(n,r)_{L_{p,q}}\hspace{1pt}$, we
choose the homology basis\/ $(\alpha_1, \beta_1) =
(\mu(A),\lambda(A))$, and\/ $(\alpha_2, \beta_2) =
(\lambda(\hat{B}_{p,q}), -\mu(\hat{B}_{p,q}))$.  The
identification with the fiber sum is now immediate.
\end{proof}

\section{Seiberg-Witten Invariants}\label{section:sw}

To prove Theorem~\ref{theorem:1}, we need to show that the family
of tori $\{T_{p,q}\}_{p\geq 2}$\/ we constructed are mutually
non-isotopic for fixed\/ $q$.  It is enough to show that the
corresponding family of 4-manifolds $\{ E(n,1)_{L_{p,q}} \}_{p\geq
2}\,$ are mutually non-diffeomorphic, since any isotopy between
the tori $\{T_{p,q}\}$ will translate into diffeomorphism between
corresponding 4-manifolds $\{ E(n,1)_{L_{p,q}} \}$.  We shall
compute their Seiberg-Witten invariants to show that $\{
E(n,1)_{L_{p,q}} \}$ are mutually non-diffeomorphic.

Let\/ $\Delta_{L_{p,q}}(x,t)$\/ denote the Alexander polynomial of
the link\/ $L_{p,q}\hspace{1pt}$, where the variable\/ $x$\/
corresponds to the axis\/ $A$.
\begin{lemma}\label{lemma:alexander}
$\Delta_{L_{p,2}} (x,t) = 1+xt^{2p-1}\hspace{1pt},$ and when $\,q
\geq 3\hspace{1pt},$
\[ \Delta_{L_{p,q}} (x,t) = 1+x^{q-1}t^{2p+q-3}+xt^2
\left( \frac{1+t^{2p-3}}{1+t} \right)
\left( \frac{1-(xt)^{q-2}}{1-xt} \right)\, . \]
\end{lemma}

\begin{proof}
The braid group on $q$ strands is generated by the elementary
braid transpositions\/ $\sigma_1,\dots ,
\sigma_{q-1}\hspace{1pt}$, where $\sigma_i$ denotes the crossing of
the $(i+1)$st strand over the $i\hspace{.5pt}$th. Note that
\[
B_{p,q} \,=\: \sigma_{q-1}\,\sigma_{q-2} \,\cdots \,
\sigma_2 \,\sigma_1^{2p-1} \, .
\]
By Theorem~1 in \cite{morton}, we have
\begin{equation}\label{eq:zero}
\Delta_{L_{p,q}}(x,t)\,=\, \det\left(I-x\, C^{(q-1)}_{q-1}
C^{(q-1)}_{q-2} \cdots\, C^{(q-1)}_2 (C^{(q-1)}_1)^{2p-1}\right)\,
,
\end{equation}
where
$C^{(q-1)}_i$ denotes the following\/ $(q-1)\times(q-1)$\/ matrix
which differs from the identity matrix $I$\/ only in the three places shown
on the
$i\hspace{.5pt}$th row.
\[ C^{(q-1)}_i  \; :=\;
\left( \begin{array}{ccccccc} 1&&&&&&\\
&\ddots&&&&&\\&&1&&&&\\
&&t&-t&1&&\\
&&&&1&&\\
&&&&&\ddots& \\
&&&&&&1
\end{array}\right)
\!\!\begin{array}{r}
\\
\\
\\
\\
\\
\\
\\ .
\\
\end{array} \]
When\/ $i=1$\/ or\/ $i=q-1\hspace{1pt}$, the matrix is truncated
appropriately to give two non-zero entries in row $i$.

The statement for\/ $q=2$\/ case is obvious. We are going to prove
the case\/ $q \geq 3$\/ by induction on\/ $q$. The main step of
this induction is proving that $ D_{p,q}= xt D_{p,q-1} $, where
$D_{p,q}:= \Delta_{L_{p,q+1}}(x,t)-\Delta_{L_{p,q}}(x,t)$. Then we
prove that $\Delta_{L_{p,3}}(x,t)= 1 + x^2t^{2p} +xt^2 \left(
\frac{1+t^{2p-3}}{1+t} \right)$ and this implies $D_{p,2}=
xt^2(1-t+t^2- \cdots - t^{2p-3})+x^2t^{2p} $ and hence $D_{p,q}=
x^{q-1}t^q (1-t+t^2- \cdots - t^{2p-3}) + x^qt^{2p+q-2}$. Finally,
by using the induction assumption and\/ $\Delta_{L_{p,q+1}} (x,t)
= \Delta_{L_{p,q}}(x,t) + D_{p,q}$ we finish the proof of the
lemma.

$\Delta_{L_{p,q+1}}(x,t)= \det (I - x \Gamma_{p,q+1})$, where\/
$\Gamma_{p,q+1} := C^{(q)}_q C^{(q)}_{q-1} \cdots C^{(q)}_2
(C^{(q)}_1)^{2p-1}$. Note that
$$ C^{(q)}_i \;=\;
\left( \begin{array}{cccc} &&&0\\
&C^{(q-1)}_i&&\vdots\\
&&&0 \\
0&\dots&0&1
\end{array}\right) $$
for $i \in \{ 1,2, \dots , q-2 \},\;$ so we must have
\vspace{-10pt}
\begin{eqnarray*}
\Gamma_{p,q+1} &=& C^{(q)}_q C^{(q)}_{q-1}
\left( \begin{array}{cccc} &&&0\\
&C^{(q-1)}_{q-1}&&\vdots\\
&&&0 \\
0&\dots&0&1
\end{array}\right)^{\hspace{-8pt}
\begin{array}{l}  -1 \end{array}}
\hspace{-5pt}\left( \begin{array}{cccc} &&&0\\
&\Gamma_{p,q}&&\vdots\\
&&&0 \\
0&\dots&0&1
\end{array}\right) \\
&=&
\left( \begin{array}{cccc} 1&&&\\
&\ddots&&\\
&&1&1 \\
&&t&0
\end{array}\right)
\left( \begin{array}{cccc} &&&0\\
&\Gamma_{p,q}&&\vdots\\
&&&0 \\
0&\dots&0&1
\end{array}\right)\!\!\begin{array}{r}
\\
\\
\\ .
\end{array}
\end{eqnarray*}
Hence it follows that
\begin{equation}\label{eq:one} \Gamma_{p,q+1} \;=\;
\left( \begin{array}{cc} &0\\
\Gamma_{p,q}&\vdots\\
&0\\
&1 \\
t \, (\Gamma_{p,q})_{( q-1\hspace{1pt} ,\,\ast\hspace{1pt} )}&0
\end{array}\right)
\end{equation}
and
$$ I-x \Gamma_{p,q+1} \;=\;
\left( \begin{array}{cc} &0\\
I-x \Gamma_{p,q}&\vdots\\
&0\\
&-x \\
-xt \, (\Gamma_{p,q})_{(q-1\hspace{1pt} ,\,\ast\hspace{1pt})} &1
\end{array}\right)
\!\!\begin{array}{r}
\\
\\
\\
\\
\\ ,
\end{array}
$$
where\/ $(\Gamma_{p,q})_{(q-1\hspace{1pt} ,\,\ast\hspace{1pt})}$\/
denotes the last
row of $\Gamma_{p,q}\hspace{1pt}$.
When we calculate the determinant of the matrix $I - x \Gamma_{p,q+1}$
by expanding along its last column we get the following equality:
\begin{equation}\label{eq:two}
\; \det (I-x\Gamma_{p,q+1} ) =
\det ( I-x\Gamma_{p,q}) - (-x) t \left[\det (I-x\Gamma_{p,q}) -
\det (I-x\Gamma_{p,q-1})\right]\, .
\end{equation}

To prove the above equality, observe that all but the last row of
the minor of the matrix $I-x\Gamma_{p,q+1}$ corresponding to the
entry $-x$ in the last column are the same as the rows of
$I-x\Gamma_{p,q}\hspace{1pt}$, and the last row of the minor is
$t$\/ times the last row of $I-x\Gamma_{p,q}$ except for the last
entry. In the minor, this entry is $0$, whereas in
$I-x\Gamma_{p,q}\hspace{1pt}$ this entry is $1$ (since
Equation~(\ref{eq:one}) shows that the last diagonal entry of
$\Gamma_{p,q}$ is $0$)\footnote{Strictly speaking, the last
diagonal entry of $\Gamma_{p,q}$ is $0$, only if $q \geq 4$.
Nevertheless, for $q=3$, Equation~(\ref{eq:two}) could be
confirmed by a direct calculation of \/$\Delta_{L_{p,4}}(x,t)$
along the same lines as the calculation of
\/$\Delta_{L_{p,3}}(x,t)$.}. This observation is why the
determinant of the minor corresponding to $-x$ is $t$ times the
difference between the determinant of $I-x\Gamma_{p,q}$ and the
determinant of the minor of $I-x\Gamma_{p,q}$ obtained by deleting
the last row and the last column (and this minor is nothing but
$I-x\Gamma_{p,q-1}$).

$ D_{p,q}= xt D_{p,q-1}$\/ now follows trivially from
Equation~(\ref{eq:two}). The next step is to calculate
\/$\Delta_{L_{p,3}}(x,t)$. By Equation~(\ref{eq:zero}), we have
$$\Delta_{L_{p,3}}(x,t) \:=\: \det \left[ \left( \begin{array}{cc}
  1 & 0 \\
  0 & 1
\end{array} \right) -x \left(
\begin{array}{cc}
  1 & 0 \\
  t & -t
\end{array} \right) \left(
\begin{array}{cc}
  -t & 1 \\
  0 & 1
\end{array} \right)^{2p-1} \right]\!\!\begin{array}{l} \\ . \end{array}
$$
An easy induction argument shows that
$$\left( \begin{array}{cc}
  -t & 1 \\
  0 & 1
\end{array} \right)^{2p-1} = \;\left( \begin{array}{cc}
  -t^{2p-1} & \frac{1+t^{2p-1}}{1+t} \\
  0 & 1
\end{array} \right)$$
and then an explicit calculation gives $\Delta_{L_{p,3}}(x,t) = 1
+ x^2t^{2p} + xt^2 \left( \frac{1+t^{2p-3}}{1+t} \right)$. As we
discussed before, to finish the proof one only needs to check that
\begin{eqnarray*}
1\!\!&+&{}\!\! x^{q-1}t^{2p+q-3}
\;+\; xt^2\left( \frac{1+t^{2p-3}}{1+t} \right)
\left( \frac{1-(xt)^{q-2}}{1-xt} \right) \\[5pt]
 &+& {}\!\! x^{q-1}t^q
(1-t+t^2- \cdots - t^{2p-3}) \;+\; x^qt^{2p+q-2}
\end{eqnarray*}
\[ =\; 1\;+\; x^qt^{2p+q-2} \;+\; xt^2\left( \frac{1+t^{2p-3}}{1+t} \right)
\left( \frac{1-(xt)^{q-1}}{1-xt} \right)\!\!\begin{array}{l} \\ .
\end{array} \hfill\qedhere \]
\end{proof}

Recall that the Seiberg-Witten invariant\/
$\overline{SW}_{\!\!X}$\/ of a 4-manifold $X$\/
(with\/ $b_2^+(X)>1$)
can be thought of
as an element of the group ring of $H_2(X;\zz)$, i.e.\/
$\overline{SW}_{\!\!X} \in \zz [ H_2( X ; \zz) ]$.
If we write\/
$\overline{SW}_{\!\!X} = \sum_g a_g g \hspace{1pt}$, then we
say that\/ $g\in H_2( X ; \zz)$\/ is a Seiberg-Witten \emph{basic class}\/
of $X$\/ if\/ $a_g\neq 0$.

\begin{lemma}\label{lemma:sw}
Let\/ $\iota:[S^1\times(S^3\setminus\nu L_{p,q})]\rightarrow
E(n,r)_{L_{p,q}}$\/ be the inclusion map.  Let\/
$\xi:=\iota_{\ast}[S^1\times\mu(A)],$
$\tau:=\iota_{\ast}[S^1\times\mu(\hat{B}_{p,q})]\in
H_2(E(n,r)_{L_{p,q}} ;\zz).$  If\/ $p\geq 2$ and\/ $q\geq 3$, then
the Seiberg-Witten invariant of\/ $E(n,r)_{L_{p,q}}$ is
\begin{eqnarray*}
\overline{SW}_{\!E(n,r)_{L_{p,q}}} \hspace{-5pt} &=& \hspace{-2pt}
(\xi^{-1}-\xi)^{n-1} \,(\xi^{-q}-\xi^q)^{r-1} \,\Big[\:
\xi^{-(q-1)}\tau^{-(2p+q-3)} +\: \xi^{q-1}\tau^{2p+q-3}
\\
&+& \hspace{-5pt} \Big(\tau^{-(2p-4)}
\sum_{i=0}^{2p-4}(-\tau^2)^{i} \Big)\, \Big( \,\sum_{j=0}^{q-3}
(\xi\tau)^{-(q-3)+2j} \Big)\:\Big]\, .
\end{eqnarray*}
When\/ $p=1$\/ and\/ $q\geq2$, we have
\[
\overline{SW}_{\!E(n,r)_{L_{1,q}}} = \; (\xi^{-1}-\xi)^{n-1}
\,(\xi^{-q}-\xi^q)^{r-1}\, \Big[ \sum_{j=-1}^{q-2}
(\xi\tau)^{-(q-3)+2j} \, \Big]\, .
\]
When\/ $q=2$\/ and\/ $p\geq1$, we have
\[
\overline{SW}_{\!E(n,r)_{L_{p,2}}} =\; (\xi^{-1}-\xi)^{n-1}
\,(\xi^{-2}-\xi^2)^{r-1}\, \left[\, \xi^{-1}\tau^{-(2p-1)} +\:
\xi\tau^{2p-1} \right] .
\]
\end{lemma}

\begin{proof}
Recall from \cite{fs:blowdown} that\/ $\overline{SW}_{\!\!E(n)} =
([F]^{-1}-[F])^{n-2}$.
From the formulas in \cite{doug:pft3} and
\cite{Taubes:T^3}, we know that
\[
\overline{SW}_{\!\!E(n)\setminus\nu
F}\:=\:([F]^{-1}-[F])\cdot\overline{SW}_{\!\!E(n)} =\:
([F]^{-1}-[F])^{n-1}
\]
From the formulas in \cite{fs:knots} and \cite{Taubes:T^3}, we conclude that
\[
\overline{SW}_{\!E(n,r)_{L_{p,q}}} =\;
\overline{SW}_{\!\!E(n)\setminus\nu F_1}\,\cdot\,
\overline{SW}_{\!\!E(r)\setminus\nu F_2}\,\cdot\,
\Delta_{L_{p,q}}^{sym} (\xi^2 , \tau^2)\, .
\]
We also need to identify\/ $[F_1]$\/ with\/ $\xi$\/ and\/ $[F_2]$\/
with\/ $\xi^q$. This last identification
is really necessary since we specifically chose to have
$\,\alpha_2 = \lambda(\hat{B}_{p,q})=q\mu(A) =q\alpha_1\in
H_1(S^3\setminus\nu L_{p,q})\,$
in the proof of Lemma~\ref{lemma:main}. Hence we must
have $[S^1\times\alpha_2]=q[S^1\times\alpha_1]\in
H_2(S^1\times(S^3\setminus \nu L_{p,q}))$.  Finally we note that
$$\Delta_{L_{p,q}}^{sym}(x,t) \:=\:
x^{-\frac{q-1}{2}}t^{-\frac{2p+q-3}{2}} \Delta_{L_{p,q}}(x,t)\,
,$$ where\/ $\Delta_{L_{p,q}}(x,t)$\/ is presented as in
Lemma~\ref{lemma:alexander}.  The rest of the proof is an easy
exercise which we shall leave to the reader.
\end{proof}

\begin{corollary}\label{cor:non-diffeo}
For fixed\/ $q\geq 3$\/ and\/ $n\geq 1$, the $4$-manifolds $\,\{
E(n,1)_{L_{p,q}} \}_{p\geq 2}\,$ are mutually non-diffeomorphic.
\end{corollary}

\begin{proof}
We note that the total number of the $SW$\/ basic classes are
invariant under any diffeomorphism.  As a consequence of Lemma~\ref{lemma:sw},
the total number of basic classes of our $4$-manifolds depend on\/
$p$. In fact, for\/ $q\geq 3$\/ and\/ $p\geq 2$, the total number
of basic classes of $E(n,1)_{L_{p,q}}$ is $(2n+2q-6)p+(qn-4n-4q+12)$.
One could check the validity of this formula by  the following
elementary argument which occupies the rest of the proof.

We first recall that the
homology classes\/ $\xi$\/ and\/ $\tau$\/ are linearly independent in
$H_2(E(n,r)_{L_{p,q}})$\/ by Proposition~3.2 in
\cite{McMullen-Taubes}.
Let\/ $N_{n,p,q}$\/ be the number of basic classes of\/ $E(n,1)_{L_{p,q}}$.
Then $N_{n,p,q}$ is the number of nonzero coefficients of
\begin{eqnarray*}
\overline{SW}_{\!E(n,1)_{L_{p,q}}} \hspace{-5pt} &=& \hspace{-2pt}
\Big[\sum_{k=0}^{n-1} (-1)^k \left(\hspace{-5pt}\begin{array}{c} n-1 \\
k\end{array}\hspace{-5pt}\right)
 \xi^{n-2k-1}\Big] \,\Big[\:
\xi^{-(q-1)}\tau^{-(2p+q-3)} +\: \xi^{q-1}\tau^{2p+q-3}
\\
&+& \hspace{-5pt} \Big(\tau^{-(2p-4)}
\sum_{i=0}^{2p-4}(-\tau^2)^{i} \Big)\, \Big( \,\sum_{j=0}^{q-3}
(\xi\tau)^{-(q-3)+2j} \Big)\:\Big]
\\
&=& \Big[ \sum_{k=0}^{n-1} (-1)^k
\left(\hspace{-5pt}\begin{array}{c} n-1 \\
k\end{array}\hspace{-5pt}\right) ( \xi^{n-2k-q}\tau^{3-2p-q}
+\xi^{n+q-2k-2}\tau^{2p+q-3} ) \Big] \\
&+& \Big[ \sum_{k=0}^{n-1} \sum_{i=0}^{2p-4}\sum_{j=0}^{q-3}
(-1)^{i+k} \left(\hspace{-5pt}\begin{array}{c} n-1 \\
k\end{array}\hspace{-5pt}\right) \xi^{n+2j+2-2k-q}
\tau^{2i+2j+7-2p-q} \Big] \, .
\end{eqnarray*}
Note that, if\/ $n+2j+2-2k-q=n+2j'+2-2k'-q$\/ and\/
$2i+2j+7-2p-q=2i'+2j'+7-2p-q$, then\/ $i+k = i'+k'$\/ hence\/
$N_{n,p,q}- 2n$ is equal to the number of elements in the set
$$\Lambda \: :=\: \{\,(j-k,i+j)\;\big| \;  0\leq i \leq 2p-4
\, ,\; 0\leq j \leq q-3 \, , \; 0\leq k\leq n-1  \,\}\, .$$
On the other hand, the number of elements in $\Lambda$ is
$$\#\Lambda \:=\: \big((2p-4)+1\big)\big((q-3)+1\big)\big((n-1)+1\big)-
(2p-4)(q-3)(n-1).$$
To see this, consider the linear map\/ $\rho:\rr^3\rightarrow \rr^2$ given by
$\rho(x,y,z)=(y-z,x+y)$.  The kernel of $\rho$ is the line generated by
the vector
$(-1,1,1)$.  Hence  $\#\Lambda$\/ is equal to the number
of lines in $\rr^3$ that are parallel to the vector $(-1,1,1)$ and
meet an integer point in the rectangular parallelepiped
\[
\mathcal{P} := \:\{\,(x,y,z)\; \big| \;  0\leq x \leq 2p-4
\, ,\; 0\leq y \leq q-3 \, , \; 0\leq z \leq n-1  \,\}\, .
\]
We  see easily that the set of such lines is parameterized by the
integer points on
the three sides of $\mathcal{P}$\/ lying on the $xy$-plane,
$xz$-plane, and  the plane\/ $x=2p-4$.
The number of integer points on these three sides of
$\mathcal{P}$\/ is equal to the
total number of integer points in $\mathcal{P}$\/ minus the number of integer
points in a smaller parallelepiped of dimensions one less.
See Figure~\ref{fig:cube}.

\begin{figure}[!ht]
\begin{center}
\includegraphics[scale=0.5]{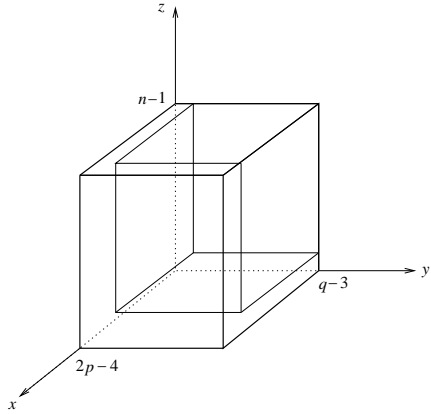}
\end{center}
\caption{The rectangular parallelepiped $\mathcal{P}$\/
in $\rr^3$}\label{fig:cube}
\end{figure}

\noindent
This proves the formula for $\#\Lambda$.
Now an easy calculation gives
\[ N_{n,p,q}\,=\: (2n+2q-6)p+(qn-4n-4q+12).  \qedhere \]
\end{proof}

It follows that the symplectic tori\/ $\{T_{p,q}\}_{p\geq 2}$\/ are
non-isotopic in $E(n)$, which completes the proof of
Theorem~\ref{theorem:1} for\/ $q\geq3$. We should note that
since the number of basic classes of\/ $E(n,1)_{L_{1,q}}$ is\/
$qn$, which is the same as the number of basic classes of
$E(n,1)_{L_{2,q}}$, the above argument cannot be used to
conclude that no torus in\/ $\{T_{p,q}\}_{p\geq 2}$\/ is isotopic
to $T_{1,q}\hspace{1pt}$. It is also insufficient to show that the
$4$-manifolds
$E(n,1)_{L_{p,2}}$\/ are mutually non-diffeomorphic,
since the number of their basic classes are independent of $p$.
Nevertheless, it is possible to prove the following corollary by comparing the
divisibilities of certain basic classes of\/ $E(n,r)_{L_{p,2}}$.

\begin{corollary}\label{cor:q=2}
There are infinitely many mutually non-isotopic tori in the set of symplectic
tori\/ $\{T_{p,2}\}_{p\geq 2}$\/ in\/ $E(n)$.
\end{corollary}

\begin{proof}
First note that by an easy Mayer-Vietoris argument
we can extend the set\/ $\{ \xi, \tau\}$\/ to an integral basis of
$H_2(E(n,r)_{L_{p,q}} ;\zz)$
(see \cite{gs} pp.$\:$72--74\/ for the description of a basis of
$H_2([E(n)\setminus\nu F]\hspace{1pt};\zz)\,$).
To distinguish\/ $T_{p_1,2}$\/ and\/ $T_{p_2,2}$\/ in\/ $E(n)$\/ we
could choose a suitable $r$ and compare the divisibilities of basic
classes with Seiberg-Witten invariant\/ $\pm 1$\/ for\/ $E(n,r)_{L_{p_i,2}}$.
When $n$ is odd and\/ $p_1>p_2>\frac{n-1}{2}$, a comparison of
the divisibilities of basic classes for\/
$E\left( n,\frac{2p_1+1-n}{2}\right)_{L_{p_i,2}}$
shows that these $4-$manifolds are not diffeomorphic and hence\/
$T_{p_1,2}$\/ and\/ $T_{p_2,2}$\/ are non-isotopic in\/ $E(n)$.
In case $n$ is even $r$ could be chosen as\/ $2p_1-\frac{n}{2}$\/ to
distinguish\/ $T_{p_1,2}$\/ and\/ $T_{p_2,2}$\/ in\/ $E(n)$\/ for
sufficiently large\/ $p_2 < p_1$.
\end{proof}

\section{Generalization to Other 4-Manifolds}\label{section:generalization}

It is not hard to show that our family\/ $\{ T_{p,q} \}$\/ embed
non-isotopically into a more general class of 4-manifolds than\/ $E(n)$.
One way to generalize Theorem~\ref{theorem:1} is to use the formulas in
\cite{fs:knots} and show that there is an infinite family of non-isotopic
symplectic tori in\/ $q[F]$\/ for\/ $q\geq3$ provided that $F$ is a
$c$-embedded symplectic torus of self-intersection $0$ in a symplectic
$4$-manifold (cf.$\;$Theorem 5.2 in \cite{fs:non-isotopic}). Another way
is to utilize the product formulas for Seiberg-Witten invariants
given in \cite{doug:pft3} and \cite{Taubes:T^3} to generalize
Theorem~\ref{theorem:1} to an even larger class of $4$-manifolds.
The latter will be our approach here.
In the rest of this section let $F$\/ denote a 2-torus of
self-intersection $0$, which is a symplectic submanifold of a
symplectic 4-manifold $X$.  We shall always assume that $[F]\in
H_2(X;\zz)$ is a primitive class, i.e. $[F]\neq m \gamma$\/ for
some\/ $\gamma\in H_2(X;\zz)$\/ and\/ $|m|>1\hspace{1pt}$.  Just
as before, we may identify a tubular neighborhood $\nu F$\/ with
$[S^1\times(S^3\setminus \nu A)]$, and we immediately obtain an
embedding of the tori\/ $\{ T_{p,q} \}$\/ into $\nu F \subset X$.
It is easy to check that Lemma~\ref{lemma:symplectic} continues to
hold in this more general setting, i.e. $T_{p,q}\subset X$\/ is a
symplectic submanifold representing the homology class $q[F]$.

To show that the family\/ $\{ T_{p,q} \}$\/ are mutually
non-isotopic in $X$, we need to again show that the corresponding
4-manifolds,
\[
X\#_{\varphi_{p,q}} E(r) \: := \; [X \setminus \nu T_{p,q}]
\hspace{-20pt} \bigcup_{T_{p,q}\times\partial D^2=F_2\times
\partial D^2} \hspace{-20pt} [E(r)\setminus \nu F_2]
\]
\vspace{-5pt}
\[
=\;\;\; [X\setminus \nu F]\hspace{-20pt}\bigcup_{F\times\partial
D^2=(S^1\times \alpha_1)\times\beta_1}\hspace{-20pt}
[S^1\times(S^3\setminus \nu
L_{p,q})]\hspace{-20pt}\bigcup_{(S^1\times
\alpha_2)\times\beta_2=F_2\times\partial D^2}\hspace{-20pt}
[E(r)\setminus\nu F_2]\; ,
\]
are mutually non-diffeomorphic.  Here we once again choose the
homology basis\/ $(\alpha_1, \beta_1) = (\mu(A),\lambda(A))$,
and\/ $(\alpha_2, \beta_2) = (\lambda(\hat{B}_{p,q}),
-\mu(\hat{B}_{p,q}))$\/ when we do link surgery.  To compute the
Seiberg-Witten invariants, we will need the following definition.

\begin{definition}\label{def:essential}
A surface\/ $F\subset X$\/ is said to be \emph{essentially
embedded}\/ if there exists a homology class\/ $\delta \in
H_2(X;\zz)$\/ such that $[F]\cdot \delta \neq 0$.
\end{definition}

Note for example that a torus fiber $F_2\subset E(r)$\/ is
essentially embedded since we can take\/ $\delta$\/ to be the
homology class of a section.

\begin{lemma}\label{lemma:SW_X}
Suppose\/ $F\subset X$\/ is an essentially embedded\/ $2$-torus of
self-intersection\/ $0$.  If\/ $H^1(X\setminus \nu F\hspace{1pt}
;\zz)=0$, then
\[
\overline{SW}_{\!\!X\#_{\varphi_{p,q}} E(r)} \,=\;
\overline{SW}_{\!\! X} \cdot ([F]^{-1}-[F])\, ([F]^{-q}-[F]^q)^{r-1}
\,
 \Delta_{L_{p,q}}^{sym} ( [F]^2 ,
\tau^2 ) \, ,
\]
where\/ $\tau :=\iota_{\ast}[S^1\times\mu(\hat{B}_{p,q})]\in
H_2(X\#_{\varphi_{p,q}}E(r)\hspace{1pt};\zz)$\/ is defined as in
Lemma~$\ref{lemma:sw}$.
\end{lemma}

\begin{proof}
Corollary~20 in \cite{doug:pft3} gives
\[
\overline{SW}_{\!\!X\setminus\nu F} =\; \overline{SW}_{\!\!X}
\cdot ([F]^{-1}-[F]).
\]
The rest of the proof goes exactly the same way as in the proof of
Lemma~\ref{lemma:sw}.
Again we need to identify\/ $[F_2]=q[F]$\/ for the same reasons as
in the proof of Lemma~\ref{lemma:sw}.
\end{proof}

\begin{corollary}\label{cor:X_pq}
Suppose that\/ $F$\/ is an essentially embedded symplectic\/
$2$-torus in a symplectic\/ $4$-manifold\/ $X$\/ with\/
$b_2^+(X)>1$.  Also assume that\/ $[F]\in H_2(X;\zz)$\/ is
primitive,\/ $[F]\cdot[F]=0$, and\/ $H^1(X\setminus\nu
F\hspace{1pt};\zz)=0$. Then for fixed\/ $q\geq 3$,
there are infinitely many manifolds in the set of
symplectic\/ $4$-manifolds $\,\{
X\#_{\varphi_{p,q}}E(1) \}_{p\geq 2}\,$ that are mutually
non-diffeomorphic.
\end{corollary}

\begin{proof}
When\/ $b_2^+(X)>1$, we always have\/ $\overline{SW}_{\!\!X}\neq
0$\/ by Taubes' result \cite{Taubes:symplectic}.  Note that the
homology classes $[F]$ and $\tau$ are linearly independent in
$H_2(X\#_{\varphi_{p,q}}E(1))$\/ by Proposition~3.2 in
\cite{McMullen-Taubes}.  The rest of the proof is an easier
analogue of the proof of Corollary~\ref{cor:non-diffeo}.
We formally set\/ $\tau =[F]$, and note that as $p\rightarrow \infty$,
the total number of terms in
$$ \overline{SW}_{\!\!X\#_{\varphi_{p,q}} E(1)}
\Big|_{\tau =[F]}
\,=\;
\overline{SW}_{\!\! X} \cdot ([F]^{-1}-[F])
\,
 \Delta_{L_{p,q}}^{sym} ( [F]^2 ,
[F]^2 )  $$
goes to $\infty$  as well.  This immediately implies that the
number of $SW$ basic classes of $X\#_{\varphi_{p,q}} E(1)$\/
also goes to $\infty$\/ as $p\rightarrow \infty$.
\end{proof}

Now we can readily generalize Theorem~\ref{theorem:1} to other
symplectic pairs\/ $(X,F)$\/ as follows.

\begin{theorem}\label{theorem:2}
Suppose that\/ $F$\/ is an essentially embedded symplectic\/
$2$-torus in a symplectic\/ $4$-manifold\/ $X$\/ with\/
$b_2^+(X)>1$.  Also assume that\/ $[F]\in H_2(X;\zz)$\/ is
primitive,\/ $[F]\cdot[F]=0$, and\/ $H^1(X\setminus\nu
F\hspace{1pt};\zz)=0$. Then for any integer $q\geq3$, there exists
an infinite family of mutually non-isotopic symplectic tori
representing the homology class $q[F]$. \hfill \qed
\end{theorem}

Note that, Corollary~\ref{cor:q=2} cannot be easily generalized, since the
divisibilities of basic classes of $X\#_{\varphi_{p,q}}E(r)$ depend
heavily on $\overline{SW}_{\!\!X}$ and $H_2(X;\zz)$.

\begin{remark}
One must take care and define
$\overline{SW}_{\!\!X}:=\overline{SW}_{\!\!X,
F}^{\hspace{1pt}\pm}$\/ when $b_2^+(X)=1$ (see \cite{fs:knots} and
\cite{doug:pft3}). When $b_2^+(X)=1$, it is not automatic that\/
$\overline{SW}_{\!\!X}$ is a finite sum and
$\overline{SW}_{\!\!X}\neq 0$\/ for a symplectic $X$. If indeed
$\overline{SW}_{\!\!X}\neq 0$\/ and is a finite sum, then
Theorem~\ref{theorem:2} will still be valid for such pair $(X,F)$.
However if\/ $\overline{SW}_{\!\!X}=0$\/ or is an infinite sum,
then there seems to be no method currently available to check
whether the tori in our family are mutually non-isotopic in $X$.
Along this line, it has been conjectured that an infinite family
of homologous but non-isotopic symplectic tori cannot occur in
$X$, when $X$\/ is a rational ruled surface with $c_1^2(X)>0$. In
fact, Sikorav has already proved this conjecture for\/ $X=\cpk$ in
\cite{sikorav} and Siebert and Tian announced a proof for $X=S^2
\times S^2$\/ and\/ $X=\cpk \# \cpkk$\/ (see \cite{st}).
\end{remark}

\begin{remark}
We should also point out that since our proof uses the product
formula to calculate the Seiberg-Witten invariants of
$4$-manifolds obtained by gluing  along the boundary of a regular
neighborhood of a torus and since such product formulas don't give
the complete picture in case the gluing occurs along higher genus
surfaces, we were able to construct infinite families of
non-isotopic tori only. In fact, there is an ongoing research to
determine the number of symplectic representatives (up to smooth
isotopy) in the homology class of symplectic surfaces of higher
genus. For example in \cite{smith}, by using purely topological
techniques, Smith was able to prove that for every odd number\/
$g\neq3$\/ there exists a symplectic $4$-manifold $X_g$\/ that
contains an infinite family of homologous but non-isotopic
connected symplectic surfaces of genus $g$. On the other hand,
Siebert and Tian announced in \cite{tian:ias} that there is a
unique symplectic surface in the homology class of a degree $d$
complex curve in\/ $\cpk$\/ provided that\/ $d\leq17$\/ and
similarly in the homology class of a complex curve of bidegree\/
$(m,n)$\/ in a Hirzebruch surface provided that\/ $m\leq7$.
\end{remark}

\smallskip
\subsection*{Acknowledgments}
The second author would like to thank Ronald Fintushel and Stefano
Vidussi for helpful discussions.
Some computations in Section~\ref{section:sw}
were verified with the aid of {\sl
Maple}$\hspace{1pt}^{\circledR}$   Version 8.

\end{document}